\documentclass[twocolumn,showpacs,preprintnumbers,amsmath,amssymb,showkeys,floatfix]{revtex4}

\usepackage{graphicx}
\usepackage{dcolumn}
\usepackage{bm}

\newcommand{\eps}{\varepsilon}
\newcommand{\ph}{\varphi}

\newcommand{\trn}{^{\rm\scriptscriptstyle T}}

\newcommand{\mR}{\mathbb R}

\newcommand{\inr}{\!\in \mR}

\DeclareMathOperator{\sign}{sign}

\newcommand{\be}{\begin{equation}}
\newcommand{\ee}{\end{equation}}

\renewcommand{\dfrac}{\frac}

\begin{document}

\preprint{APS}

\title{Controlled Synchronization Under Information Constraints}

\author{Alexander L. Fradkov, Boris Andrievsky}
  \email{fradkov@mail.ru; bandri@yandex.ru}%
\affiliation{Institute for Problems of Mechanical Engineering,
Russian Academy of Sciences, \\ 61, Bolshoy V.O. Av., 199178,
Saint Petersburg, Russia %
}%
\author{Robin J. Evans}
\email{r.evans@ee.unimelb.edu.au.}%
\affiliation{National ICT Australia and\\%
Department of Electrical and Electronic Engineering,\\%
University of Melbourne, Victoria, 3010, Australia%
}%
\date{\today}

\begin{abstract}
The class of controlled synchronization systems under
information constraints imposed by limited information capacity of
the coupling channel is analyzed. It is shown that the framework proposed in A. L. Fradkov, B. Andrievsky, R. J. Evans, Physical
Review E 73, 066209 (2006) is suitable not only for observer-based synchronization but also for controlled master-slave synchronization via communication channel with limited information capacity.
A simple first order coder-decoder scheme is proposed and
 a theoretical analysis for 
multi-dimensional master-slave systems represented in the  Lurie
form (linear part plus nonlinearity depending only on measurable
outputs) is provided. An output feedback control law is proposed based on the Passification theorem.
It is shown that the upper bound of the limit
synchronization error  is proportional to the upper bound of the
transmission error. As a consequence, both upper and lower bounds
of limit synchronization error are proportional to the maximum
rate of the coupling signal and inversely proportional to the
information transmission rate (channel capacity).  The results are
applied to controlled synchronization of two chaotic Chua systems
coupled via a controller and a channel with limited capacity.
\end{abstract}


\maketitle

\clearpage
\section{Introduction}
Synchronization of nonlinear systems, particularly chaotic systems has attracted the attention of many researchers for several decades \cite{Blekhman88,Pikovsky01}. During recent years interest in controlled synchronization has increased, partly driven by a growing interest in the application of control theory methods in physics \cite{PecoraCarroll90,BFNP97,FP98,Fradkov07}. The design of links interconnecting parts of complex systems to enable synchronization was studied in \cite{Belykh05,Zhan07,Cruz07}. These and related papers explore the possibility of modification of complex dynamical system behavior by means of feedback action. 

Modifying the behavior of complex interconnected systems and networks has attracted considerable interest. The available results significantly depend on models of interconnection between nodes. In some works  the interconnections are modeled as delay elements. However, the spatial separation between nodes  means that modeling connections via communication channels with limited capacity is more realistic. Hence analysis of the overall system should include both dynamical and information aspects. 

Recently the limitations of control under constraints imposed by a finite capacity information channel have been investigated in detail in the control theory literature, see \cite{WongBrockett_AC97e,NairEvans_CDC02,NairEvans_Aut03,%
NairEvans_SIAM04,NairEvans_AC04,BazziMitter_InTh05,NairFagnani07} and references
and the references therein. It has been shown that stabilization under information constraints is possible if and only if the capacity of the information channel exceeds the entropy production of the system at the equilibrium \cite{NairEvans_Aut03,NairEvans_SIAM04,NairEvans_AC04}. In \cite{Lloyd00,Lloyd04} a general statement was proposed, claiming that the difference between the entropies of the open loop and the closed loop systems cannot exceed the information introduced by the controller, including the transmission rate of the information channel. However, results of the previous works on control systems analysis under information constraints do not apply to synchronization systems since in a synchronization problem trajectories in the phase space converge to a set (a manifold) rather than to a point, i.e. the problem cannot be reduced to simple stabilization.

The first results on synchronization under information constraints were presented in \cite{FradkovAndrievskyEvans_PRE06},
where  the so called observer-based synchronization scheme \cite{PecoraCarroll90,FNM00} was considered. In this paper we extend this work and analyze a controlled synchronization scheme for two nonlinear systems. A major difficulty with the controlled synchronization problem arises because the coupling is implemented in a restricted manner via the control signal. Key tools  used to solve the problem are quadratic Lyapunov functions and passification methods borrowed from control theory. 
To reduce technicalities we restrict our analysis by Lurie systems (linear part plus nonlinearity depending only on measurable
outputs).

 \section{Description of controlled synchronization scheme}

Consider two identical dynamical systems modeled in Lurie form (i.e. the right hand sides are split into a linear part and a nonlinear part which depends only on the measurable outputs). Let one of the systems be controlled by a scalar control function $u(t)$ whose action is restricted by a vector of control efficiencies $B$. The controlled system model is as follows:
\begin{align} 
&\dot x(t)=Ax(t)+B\ph (y_1),~y_1(t)=Cx(t), \label{2}
\\
&\dot z(t)=Az(t)+B\ph (y_2)+Bu,~y_2(t)=Cz(t), \label{1}
\end{align} 
where $x(t)$, $z(t)$ are $n$-dimensional (column) vectors of state
variables; $y_1(t)$, $y_2(t)$ are scalar output variables; $A$ is 
an $(n\times n)$-matrix; $B$ is $n\times 1$ (column) matrix; $C$ is  an $1\times n$ (row) matrix, $\ph(y)$ is a continuous nonlinearity, acting in the span of control; vectors $\dot x$, $\dot z $ stand for time-derivatives of $x(t)$, $z(t)$ correspondingly. System (\ref{2}) is called {\it master} ({\it leader}) {\it system}, while the controlled system (\ref{1}) is {\it slave system} ({\it follower}). Our goal is to evaluate limitations imposed on the synchronization precision by limited the transmission rate between the systems. The intermediate problem is to find a control function ${\cal U}(\cdot)$ depending on measurable variables such that the  synchronization error $e(t),$ where $e(t)=z(t)-x(t)$ becomes small as $t$ becomes large. We are also interested in the value of output synchronization error $\eps(t)=y_2(t)$ $-y_1(t)=Ce(t)$.

\begin{figure}[htpb]
\centering
\includegraphics[width=85mm]{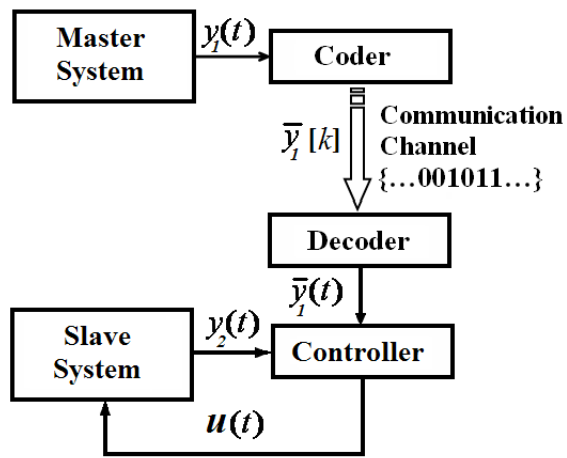}
\caption{Block-diagram for master--slave controlled synchronization (master system output $y_1$ is transmitted over the channel).}
\label{Stuc_contra}
\end{figure}

A key difficulty arises because the output of the master system is not available directly but only through  a communication channel with limited capacity. This means that the signal $y_1(t)$ must be coded at the transmitter side and codewords are then transmitted with only a finite number of symbols per second thus introducing some error. We assume that the observed signal $y_1(t)$ is coded with symbols from a finite alphabet at discrete sampling time instants $t_k=k T_s$, $k=0,1,2,\dots$, where $T_s$ is the sampling time. Let the coded symbol $\bar {y}_1 [k]=\bar {y}_1 (t_k)$ be
transmitted over a digital communication channel with a finite capacity. To simplify the analysis,
we assume that the observations are not corrupted by observation noise; transmission delay 
and transmission channel distortions may be neglected. Therefore, the discrete communication channel with  sampling period $T_s$ is considered, but it is assumed that the coded symbols are available at the receiver side at the same sampling instant $t_k=kT_s$, as they are generated by the coder.
Assume that {\it zero-order extrapolation} is used to convert the digital sequence $\bar {y}_1[k]$ to the continuous-time input of the response  system $\bar {y}_1(t)$, namely, that $\bar {y}_1(t)=\bar{y}_1[k]$ as $ kT_s \le t< (k+1)T_s $. Then the {\it transmission error} is defined as follows: 
\begin{align}
\delta_y(t)= y_1(t)- \bar {y}_1 (t).
\label{dy}
\end{align}
On the receiver side the signal is decoded introducing additional error and the controller can use only the signal $\bar y_1(t)=y_1(t)+\delta_y(t)$ instead of $y(t)$. A block diagram of the system is shown in Fig.~\ref{Stuc_contra}. 

We restrict consideration to  simple control functions in the form of static linear feedback
\begin{align} 
u(t)=-K\eps(t), \label{3}
\end{align} 
where $\eps(t)=y_2(t)$ $-y_1(t)$ denotes an output synchronization error, $K$ is a scalar controller gain. The problem of finding static output feedback 
even for linear systems is one of the classical problems of control theory. Although substantial effort has been devoted to its solution and various necessary and sufficient conditions for stabilizability by static output feedback have been obtained, most existing conditions are not testable practically \cite{Syrmos97,Stefanovski06}. Since we are dealing with a nonlinear problem further complicated by information constraints, we restrict our attention to sufficient conditions for solvability of the problem and evaluate upper bounds for synchronization error. 
To this end we introduce an upper bound on the limit synchronization error $Q=\sup \overline{\lim\limits_{t\to\infty}} \|e(t) \|$, where the supremum is taken over all admissible transmission errors.
In the next two sections the coding and decoding procedures are described and a bound on admissible transmission errors $\delta_y(t)$ is evaluated.

\section{Coding procedures }\label{Sec:coding}

In the paper \cite{FradkovAndrievskyEvans_PRE06} the properties of observer-based synchronization for Lurie systems over a limited data rate communication channel with one-step memory time-varying coder is studied.  It is shown that an upper bound on the limit synchronization error  is proportional to a certain upper bound on the transmission error. Under the assumption that a sampling time may be properly chosen, optimality of the binary coding in the sense of demanded transmission rate is established, and the relationship between synchronization accuracy and an optimal sampling time is found. On the basis of these results, the present paper deals with a binary coding procedure.

At first, we introduce the memoryless (static) binary coder to be a discretized map $q: \mR\to\mR$  as 
\begin{align}
\label{qu1}
q (y,M)=M\sign(y), 
\end{align}
where $\sign(\cdot)$ is the {\it signum} function: $\sign(y)=1$, if $y\ge 0$, $\sign (y)=-1$, if $y<0$; parameter $M$ may be referred to as a {\it coder range}. Evidently, $|y-q(y,M)|\le M$ for all $y$ such that $y: |y|\le 2M$. Notice that for a binary coder each codeword symbol contains one bit of information. The discretized output of the considered coder is given as $\bar y=q(y,M)$ and we assume that the coder and decoder make decisions based on the same information. 

The static coder \eqref{qu1} is a subset of the class of time-varying coders with memory, see e.g. \cite{NairEvans_Aut03,BrockettLiberzon_AC00,%
Liberzon03,TatikondaMitter_AC04,FradkovAndrievskyEvans_PRE06,NairFagnani07}. Two underlying ideas are used for these kinds of encoder:
\begin{itemize}
\item[-- ] reducing the coder range $M$ to cover the same area around the predicted value for the $(k+1)$th observation $y[k+1]$, $ y[k+1]\in {\cal Y}[k+1]$. This means that the quantizer range $M$ is updated with time and a time-varying quantizer (with different values of $M$ for each instant, $M=M[k]$) is used. Using such a ``zooming'' strategy it is possible to increase coder accuracy in the steady-state mode, and, at the same time, to prevent coder saturation at the beginning of the process; 
\item[-- ] introducing memory into the coder, which makes it possible to predict the $(k+1)$th observation $y[k+1]$ with some accuracy and, therefore, to transmit over the channel only the encrypted innovation signal.
\end{itemize}

In this paper a {\it first-order coder} is considered, where the predicted value $y[k+1]$ is taken equal to $y[k]$ \cite{SavkinPetersen03,TatikondaMitter_AC04,FradkovAndrievskyEvans_PRE06}. In order to describe it, introduce the sequence of {\it central numbers} $c[k]$, $k=0,1,2,\dots$ with initial condition $c[0]=0$.
At step $k$  the coder compares the current measured output $y[k]$ with the number $c[k]$, forming the deviation signal $\partial y[k]=y[k]-c[k]$. Then this signal is discretized with a given $M=M[k]$ according to \eqref{qu1}. The output signal 
\begin{align}
\label{dybar}
\bar{\partial} y[k]=q(\partial y[k],M[k])
\end{align}
 is represented as an $R$-bit information symbol from the coding alphabet
and transmitted over the communication 
channel to the decoder. Then the central number $c[k+1]$ and the range parameter $M[k]$ are renewed based on the available information about the master system dynamics. 
Assuming that the master system output $y$ changes at a slow rate,
 i.e. that $y[k+1]\approx y[k]$  
we use the following update  algorithms: 
\begin{align}
\label{ck}
c[k+1]=c[k]+\bar{\partial} y[k], \quad c[0]=0,~~k=0,1,\dots ,
\end{align}
\begin{align}
\label{mk}
M[k]=(M_0- M_{\infty})\rho^k+ M_{\infty}, ~~k=0,1,\dots ,
\end{align}
where $0<\rho\le 1$ is the decay parameter, $M_{\infty}$ stands for the limit value of $M[k]$. The initial value $M_0$ should be large enough to capture  the region of possible initial values of $y_0$. Equations \eqref{qu1}, \eqref{dybar}, \eqref{mk} describe the coder algorithm. A similar algorithm is realized by the decoder. Namely, the sequence of $M[k]$ is reproduced at the receiver node utilizing \eqref{mk}; the values of $\bar{\partial} y[k]$ are restored with given $M[k]$ from the received codeword; the central numbers $c[k]$  are found in the decoder in accordance with \eqref{ck}. Then $\bar y[k]$ is found as a sum $c[k] +\bar{\partial} y[k]$.

\section{Evaluation of synchronization error}\label{Sec:eval}
We now find a relation between the transmission rate and the achievable accuracy of the coder--decoder pair, assuming that the growth rate of $y_1(t)$ is uniformly bounded. Obviously, the exact bound ${ L_y}$ for the rate of $y(t)$ is $L_y=\sup\limits_{x\in\Omega}{|C\dot {x}|}$, where $\dot x$ is from \eqref{1}. 
To analyze the coder--decoder accuracy, evaluate the upper bound $\Delta=\sup\limits_{t}|\delta_y(t)|$ of the transmission error $\delta_y(t)=y_1(t)-\bar {y}_1(t)$. Consider the sampling interval $[t_k,t_{k+1}]$ assuming that $\partial y_1[k] \le 2M$, where $\partial y_1[k]= y_1(t_k)-c[k]$. Since $\bar{y}_1[k]=c[k]+M\sign(\partial y_1[k])$ it is clear that $|\delta_y(t_k)|\le M$. Additionally, the error may increase from $t_k$ to $t_{k+1}$ due to change of $y_1(t)$ by a value not exceeding $\sup\limits_{t_k<t<t_{k+1}}|y_1(t)-y_1(t_k)|$ $\le \int\limits_{t_k}^{t_{k+1}}{|\dot{y}_1(\tau)|d\tau}$ $ \le\int\limits_{t_k}^{t_{k+1}}{L_yd\tau}$ $=L_yT_s$. Therefore the total transmission error for each interval $[t_k,t_{k+1}]$ satisfies the inequality:
\begin{align}
|\delta_y(t)|\le M+L_yT_s
\label{delles}
\end{align}
Inequality \eqref{delles} shows that in order to meet the inequality $|\delta_y(t)|\le \Delta=2M$ for all $t$, the sampling interval $T_s$ should satisfy condition
\begin{align}
T_s<\Delta/L_y.
\label{ts}
\end{align}

Subtracting Eq. \eqref{2} from Eq. \eqref{1} and taking into account the control law  \eqref{3} we derive an equation for the synchronization error in the form
\begin{align}
\dot e(t)=\!A_Ke(t)+\!B\zeta (t) \!-\!BK\delta_y(t),
\label{10}
\end{align}
where $A_K=A-BKC$, $\zeta(t)= \ph(y_2(t)) \!-\!\ph(y_1(t)).$

 Evaluate the total guaranteed synchronization error $Q=\sup\overline{\lim\limits_{t\to\infty}}\|e(t)\|$, where $\|\cdot\|$ denotes the Euclidian norm of a vector, and the supremum is taken over all admissible transmission errors $\delta_y(t)$ not exceeding the level $\Delta$ in absolute value. The ratio $C_e=Q/\Delta$ (the relative error) can be interpreted as the norm of the transformation from the input function $\delta_y(\cdot)$ to the output function $e(\cdot)$ generated by the system \eqref{10}. Owing to the nonlinearity of the equation \eqref{10} evaluation of the norm $C_e$ is nontrivial and it even may be infinite for rapidly growing nonlinearities $\varphi (y)$. To obtain a reasonable upper bound for $C_e$ we assume that the nonlinearity is Lipschitz continuous along all the trajectories of the drive system \eqref{1}.  More precisely, we assume existence of a positive number $L_\ph>0$ such that
\begin{align*}
|\ph(y)-\ph(y+\delta)|\le L_\ph |\delta |
\end{align*}
for all $y=Cx$, $x\in\Omega$, where $\Omega$ is a set containing all the trajectories of the drive system \eqref{2}, starting from the set of initial conditions $\Omega_0$, $|\delta |\le\Delta$. For Lipschitz nonlinearities $\zeta(t)$ satisfies inequality $|\zeta(t)|\le L_{\ph}|\eps(t)|.$
After the change $K\to K+L_{\ph}$, the error equation \eqref{10} can be represented as 
\begin{align}
\dot{e}(t)=A_Ke(t)+B\xi(t)-B(K+ L_{\ph})\delta_y(t),
\label{11}
\end{align}
where the variable $\xi(t)= L_\ph \eps(t)+\zeta$, apparently, satisfies sector inequality  $\xi(t)\eps(t)\ge 0$ for all $t\ge 0$.

The problem is reduced to quantifying the stability properties of \eqref{11} for bounded input $\delta_y(t).$ We first analyze behavior of the system \eqref{11} for $\delta_y(t)=0.$ To this end, we find conditions for the existence of a quadratic Lyapunov function $V(e)=e\trn Pe$ and controller gain $K$ satisfying inequality  $\dot V(e)\le -\mu V(e)$ for some $\mu>0$ for $\delta_y(t)=0$ and for all $\xi$ satisfying the  quadratic inequality $\xi\eps\ge 0$. Such conditions can be derived from the Passification theorem \cite{Fradkov74,FradkovMN99}, see Appendix. Namely, such $V$ and $K$ exist if and only if the transfer function of the linear part of the system models  \eqref{2},  \eqref{1} $W(\lambda)=C(\lambda I-A)^{-1}B$ is {\it hyper-minimum phase} (HMP). Recall that the HMP property for a rational function $W(\lambda)=b(\lambda)/a(\lambda)$ where $a(\lambda)$ is a polynomial of degree $n$,  $ b(\lambda)$ is a polynomial of degree not greater than $n-1$ means that $ b(\lambda)$ is a Hurwitz polynomial of degree  $n-1$ with positive coefficients \cite{FradkovMN99}. Now consider the case $\delta_y(t)\ne 0$, assuming that the HMP condition holds and matrix $P$ and gain $K$ are chosen properly and the modified Lyapunov inequality $PA_K+A_K\trn P\le -\mu P$ is valid for some $\mu>0$. Evaluating the time derivative of function $V(e)$ along trajectories of \eqref{1},  \eqref{2} with initial conditions in $\Omega_0$, using standard quadratic inequality $|e\trn PB|\le\sqrt{V(e)}\sqrt{V(B)}$ after simple algebra we get
\begin{align*}
\dot V\le -\mu V+|e\trn PB(K+L_{\ph})\delta_y| \le-\mu V+\sqrt{V}\nu, 
\end{align*}
where $\nu=\sqrt{V(B)} (|K|+L_{\ph})\Delta$.
Since $\dot V\!<\!0$ within the set \mbox{$\sqrt{V}\!>\!\mu^{-1}\nu$,} the value of $\overline{\lim\limits_{t\to\infty}}\sup V(t)$ cannot exceed $\Delta^2\big(L_\ph+|K|\big)^2\lambda_{\max}(P)/\mu^2$. In view of positivity of  $P$, $\lambda_{\min}(P)\|e(t)\|^2\le V(t)$, where $\lambda_{\min}(P)$, $\lambda_{\max}(P)$ are minimum and maximum eigenvalues of $P$, respectively. Hence
\begin{align}
\overline{\lim\limits_{t\to\infty}}\|e(t)\|\le C_e^{+}\Delta,
\label{12}
\end{align}
where 
\begin{align}
C_e^{+}=\sqrt{\dfrac{\lambda_{\max}(P)}{\lambda_{\min}(P)}}\dfrac{L_\ph+|K|}{\mu}.\label{ce}
\end{align}

The inequality \eqref{12} shows that the total synchronization error is proportional to the upper bound on the transmission error $\Delta$. 

As was shown in the authors' previous paper \cite{FradkovAndrievskyEvans_PRE06}, a binary coder is optimal in the sense of bit-per-second rate, and the
optimal sampling time $T_s$ for this coder is 
\begin{align}
T_s= {\Delta}/{ (\beta L_y)},
\label{tsopt},
\end{align}
where $\beta\approx 1.688$.
Then the channel bit-rate $ R =1/T_s$ is as follows:
\begin{align}
 R = {\beta L_y}/{\Delta },
\label{ropt}
\end{align}
 and this bound is tight for the considered class of coders.
 Taking into account the relation \eqref{ropt} for optimal transmission rate, the synchronization error can be estimated as follows: 
\begin{align}
\overline{\lim\limits_{t\to\infty}}\|e(t)\|\le C_e^{+}\beta L_y/ R ,
\label{12a}
\end{align}
 i.e. it can be made arbitrarily small for sufficiently large transmission rate $R$. 

{\it Remark 1}. Related estimates for synchronization errors in coupled systems were obtained in several papers \cite{AfraimovichChowHale97,%
AfraimovichLin98,ChiuLinPeng00,%
BelykhBelykhNevidin03,BelykhBelykhHasler04}. However, in \cite{AfraimovichChowHale97,AfraimovichLin98,%
ChiuLinPeng00,BelykhBelykhNevidin03,%
BelykhBelykhHasler04} either existence of Lyapunov functions, i.e. stability of uncoupled systems is required, or a partial stability (stability of the synchronization manifold) is provided by a strong coupling playing the role of state feedback in the error system. In this paper only output feedback is allowed and coupling is applied through the control term $Bu$, i.e. in a restrictive manner. That is why the result holds under the additional assumption (passifiability) caused by the nature of controlled problems. Then the partial stability conditions are provided by  linear observer theory. In addition the final result \eqref{12a} is presented in terms of transmission rate, i.e. appeals to the information theory view. 

{\it Remark 2}. One can pose the following problem: evaluate upper and lower bounds for $C_e$ based on worst case inputs $\delta_y(t)$. Such a problem is similar to the energy control problem for systems with dissipation \cite{Fradkov05,Fradkov06} and $C_e$ can be interpreted as the {\it excitability index} of the system. Employing a lower bound for the excitability index for passive systems \cite{Fradkov05,Fradkov06} we conclude that if the gain vector $K$ is chosen to ensure strict passivity of the system \eqref{10} then the lower bound for $C_e$ is positive, i.e.
\begin{align}
\sup\limits_{|\delta_y(t)|\le\Delta}\overline{\lim\limits_{t\to\infty}}\|e(t)\|\ge C_e^{-}\Delta.
\label{12b}
\end{align}

Thus, for finite channel capacity the guaranteed synchronization error is not reduced to zero being of the same order of magnitude as the transmission error. 

{\it Remark 3.} It worth noting that relations \eqref{12}, \eqref{12b} give overestimates for the synchronization error, because they provide upper and lower bounds for the worst case of the input signal $\delta_y(t)$. Also the estimate of the mean square value of the synchronization error may be used. There is a significant  body of work in which the quantization error signal $\delta_y(t)$ is modelled as an extra additive white noise. This assumption, typical for the digital filtering theory, is reasonable if the quantizer resolution is high \cite{Curry70}, but it needs modification for the case of a low number of quantization levels \cite{NairFagnani07}.   The dependence 
\eqref{12a} will be used for numerical analysis in Sec.~\ref{Sec:example}.

\section{Example. Synchronization of chaotic Chua systems}\label{Sec:example}

Let us apply the above results to synchronization of two chaotic Chua systems coupled via a channel with limited capacity.

{\it Master system}. Let the master system \eqref{2} be represented as the following {\it Chua system}:
\begin{align}
\label{chuagen}
&\begin{cases}
\dot x_1=p(-x_1+\ph(y_1)+x_2),\quad t\ge 0,\cr
\dot x_2=x_1-x_2+x_3\cr
\dot x_3=-qx_2,
\end{cases}\\
&y_1(t)=x_1(t),\nonumber
\end{align}
where $y_1(t)$ is the master system output (to  be transmitted over the communication channel), $p$, $q$ are 
known parameters, $x=[x_1,x_2,x_3]\trn\inr^3$ is the state vector; $\ph(y_1)$ is a piecewise-linear function,  having the form:
\begin{align}
\ph(y)&= m_0y+m_1(|y+1|-|y-1|),
\label{phi}
\end{align}
where $m_0$, $m_1$ are given parameters.

{\it Slave system}. 
Correspondingly, the slave system equations \eqref{2} for the considered case become
\begin{align}
\label{chuares}
&\begin{cases}
\dot z_1=p\big(-z_1+\ph(y_2)+z_2+u(t)\big),\quad t\ge 0,\cr
\dot z_2=z_1-z_2+z_3\cr
\dot z_3=-qx_2,
\end{cases}\\
&y_2(t)=z_1(t),\nonumber
\end{align}
where $y_2(t)$ is the slave system output, $z=[z_1,z_2,z_3]\trn\inr^3$ is the state vector, $\ph(y_2)$ is defined by \eqref{phi}.

{\it Controller} has a form 
\begin{align} 
u(t)=-K\eps(t), \label{bar3}
\end{align} 
where $\eps(t)=y_2(t)$ $-\bar{y}_1(t)$; $\bar{y}_1(t)$ is a master system output, restored from the transmitted codeword  by the receiver at the slave system node (see Fig.~\ref{Stuc_contra}); the gain $K$ is a design parameter.

{\it Coding procedure} has a form \eqref{dybar}--\eqref{mk}. The input signal of the coder is $y_1(t)$. The reference input $\bar {y}_1(t)$ for controller \eqref{bar3} is found by holding the value of $\bar {y}_1[k]$ over the sampling interval $[kT_s,(k+1)T_s)$, $k=0,1,\dots$.

\begin{figure}[htpb]
\centering
\includegraphics[width=85mm]{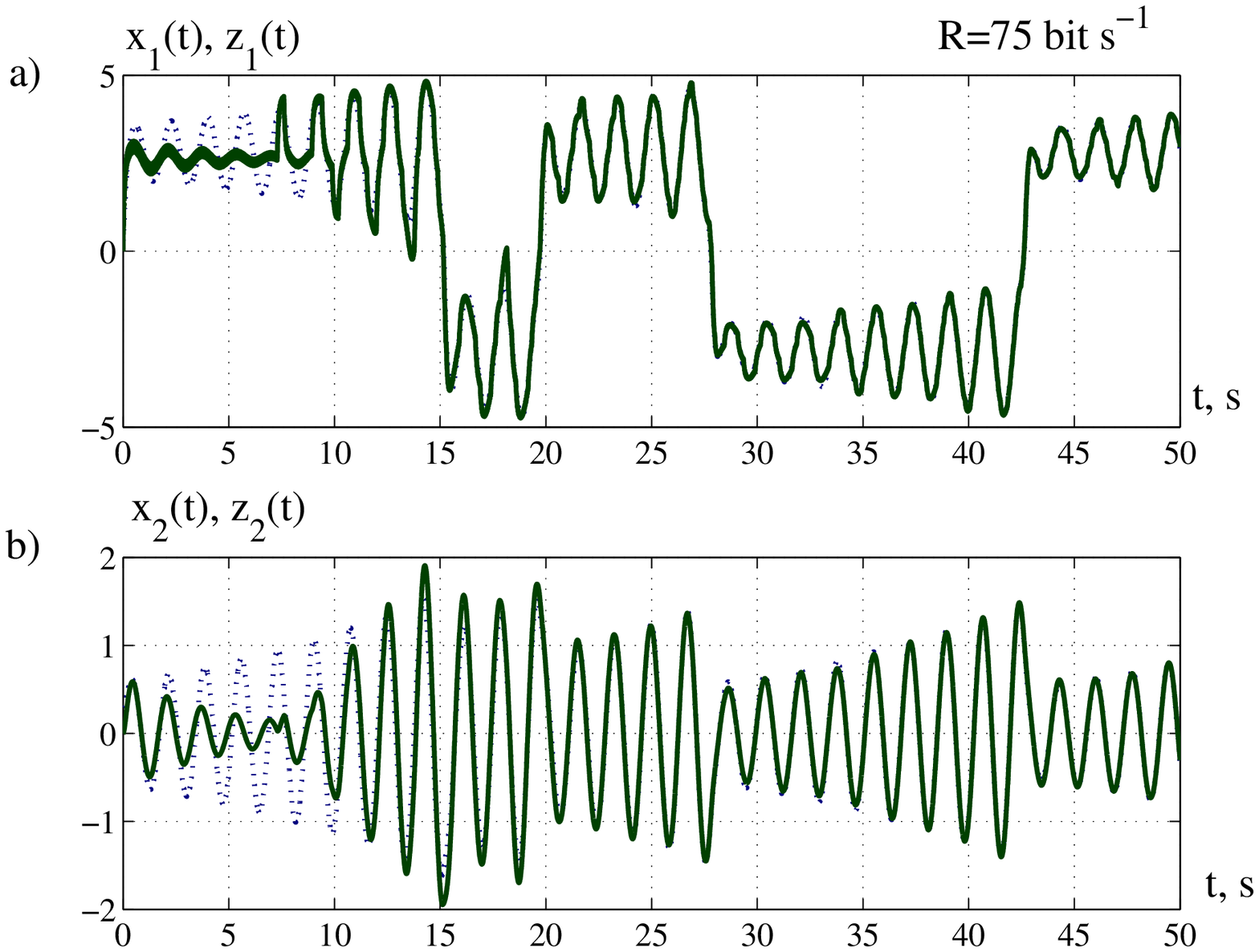}
\caption{(Color online) Time histories of the state variables of master and slave systems \eqref{chuagen}, \eqref{chuares} for $\Delta=1$ ($T_s= 13$~ms, $R=75$~bit/s): a) $x_1(t)$ (dotted
line, $z_1(t)$ (solid line); b) $x_2(t)$ (dotted
line, $z_2(t)$ (solid line).}
\label{xze050a}
\centering
\includegraphics[width=85mm]{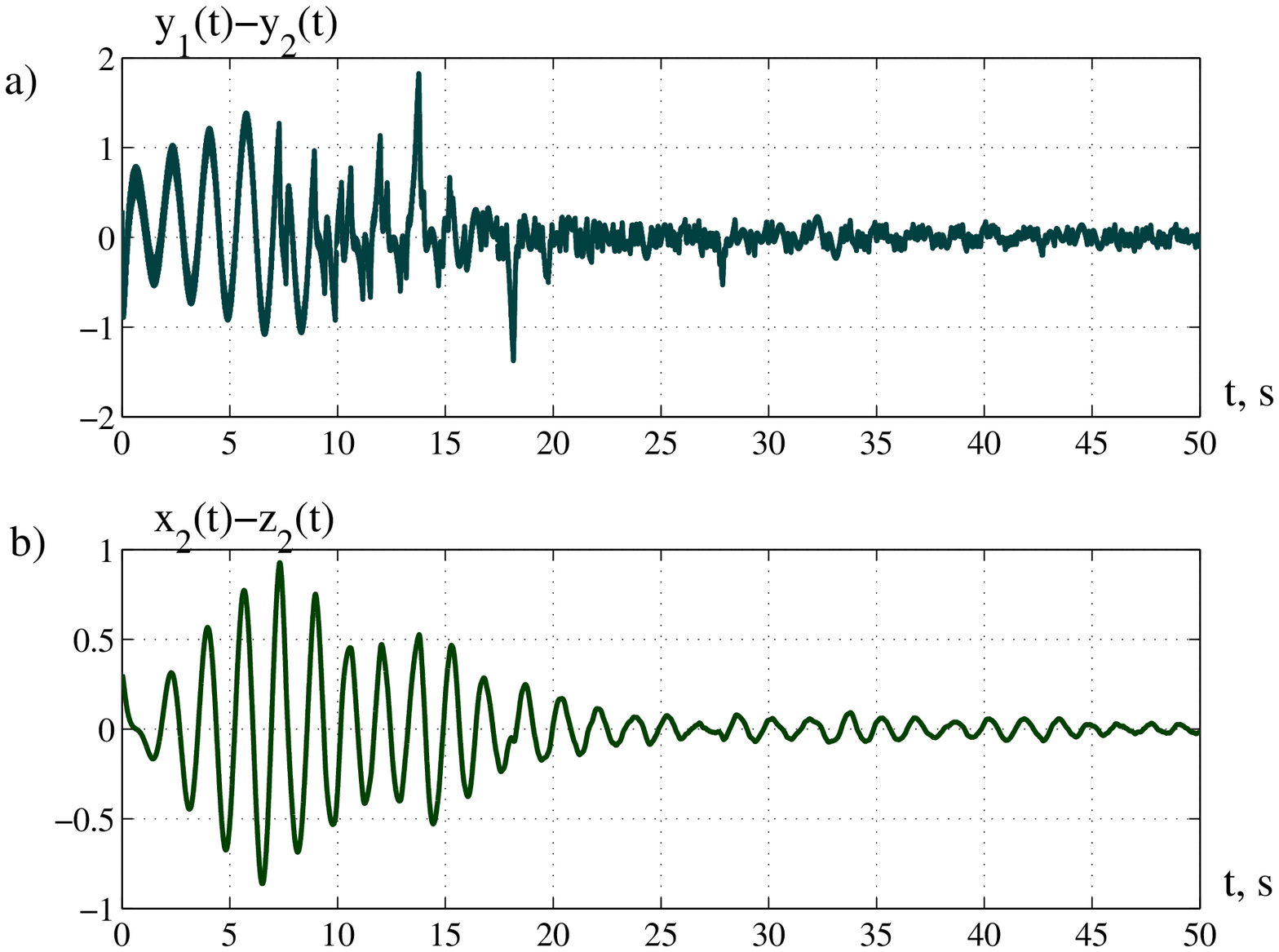}
\caption{(Color online) Synchronization error time histories: a) $\eps(t)=y_1(t)-y_2(t)$, b) $e_2(t)=x_2(t)-z_2(t)$ for $\Delta=1$ ($T_s= 13$~ms, $R=75$~bit/s).}
\label{dydx050a}
\end{figure}

\begin{figure}[htpb]
\centering
\includegraphics[width=85mm]{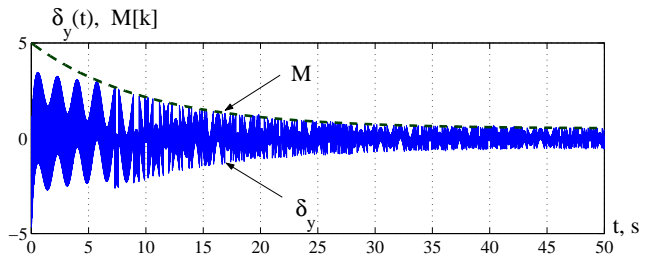}
\caption{(Color online) Time histories of transmission error $\delta_y(t)$ (solid line) and coder range parameter $M[k]$ (dashed line) for $\Delta=1$ ($T_s= 13$~ms, $R=75$~bit/s).}
\label{delyMa}
\end{figure}

\begin{figure}[htpb]
\centering
\includegraphics[width=85mm]{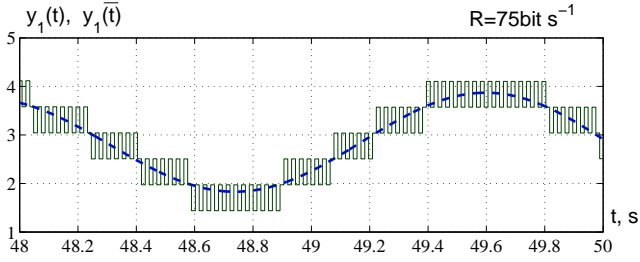}
\caption{(Color online) Time histories of $y_1(t)$ (dashed line), $\bar {y}_1(t)$ (solid line) for $R=75$~bit/s, $t\in[48, 50]$~s.}
\label{y1y1050a}
\end{figure}

\begin{figure}[htpb]
\centering
\includegraphics[width=85mm]{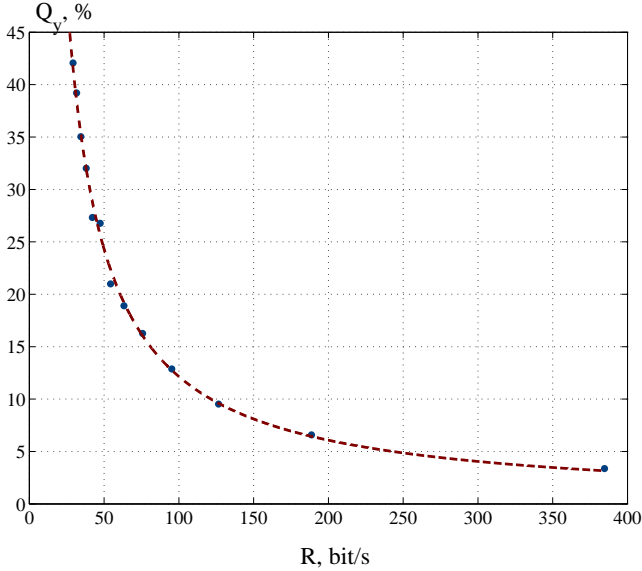}
\caption{(Color online) Normalized transmission error 
$Q_y$ vs transmission rate $ R $.}
\label{QyotRa}
\end{figure}

\begin{figure}[htpb]
\centering
\includegraphics[width=85mm]{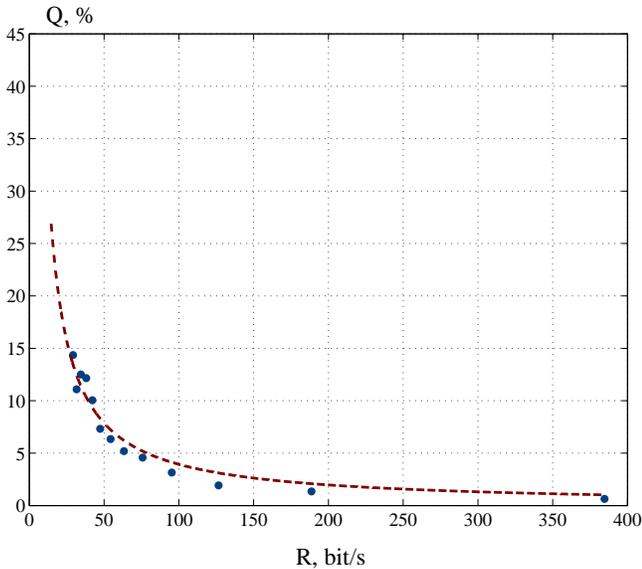}
\caption{(Color online) Normalized state synchronization error 
$Q$ vs transmission rate $ R $.}
\label{QotRa}
\end{figure}

\begin{figure}[htpb]
\centering
\includegraphics[width=85mm]{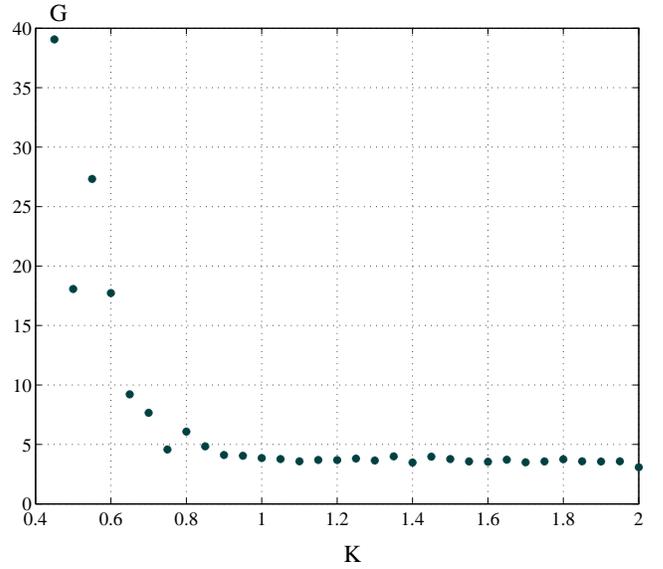}
\caption{(Color online) Overall synchronization error gain 
$G$ vs controller gain $K$.}
\label{QotKa}
\end{figure}

The following parameter values were taken for the simulation:
\begin{itemize}
\item[--] Chua system parameters: $p=10$, $q=15.6$, $m_0=0.33$, 
$m_1=0.945$; 
\item[--] the bound ${ L_y}$ for the rate of $y_1(t)$ was evaluated by numeric integration of \eqref{2} over the time interval $t\in [0,t_\text{fin}]$, $ t_\text{fin}=1000$~s, as $L_y=45$;
\item[--] parameter $\Delta$ was taken for different simulation runs as $\Delta=0.2,0.4,\dots,3.0$;
\item[--] the sample interval $T_s$ was found for each $\Delta$ from \eqref{ts};
\item[--] the following coder parameters $M_0$, $M_\infty$, $\rho$ in \eqref{mk} were taken: $M_0=5$, $ M_\infty =\Delta/2$ (different for each $\Delta$), $\rho=\exp(-0.1T_s)$;
\item[--] initial conditions for master and slave systems: $x_i =0.3$,
$z_i= 0$ ($i=1,2,3$);
\item[--] the simulation final time $t_\text{fin}=1000$~s.
\end{itemize}
The following accuracy indexes were calculated:
\par\noindent -- relative transmission error
\begin{align}
Q_y=\dfrac{\max\limits_{ 0.8t_\text{fin}\le t \le t_\text{fin}}|\delta_y(t)|}{ \max\limits_{0\le t \le t_\text{fin}}|y_1(t) |},
\label{Qy}
\end{align}
\par\noindent -- normalized state synchronization error
\begin{align}
Q=\dfrac{\max\limits_{ 0.8t_\text{fin}\le t \le t_\text{fin}}\|e(t)\|}{\max\limits_{0\le t \le t_\text{fin}}\|x(t)\|},
\label{Q}
\end{align} 
where $\delta_y(t)= y_1(t)- \bar {y}_1 (t)$, $e(t)=x(t)-z(t)$. Expressions \eqref{Qy}, \eqref{Q}  characterize transmission error and synchronization error near a steady-state mode.  

Results of the system examination for $K=1.0$ are reflected in Figs.~ \ref{xze050a}--\ref{QotRa}. 
Figure \ref{xze050a} shows time histories of the state variables of the master and slave systems \eqref{chuagen}, \eqref{chuares} $x_1(t)$, $z_1(t)$ (plot {\it a}) and $x_2(t)$, $z_2(t)$ (plot {\it b}) for $\Delta=1$ ($T_s= 13$~ms, $R=75$~bit/s). The corresponding synchronization errors $e_1(t)$ and $e_2(t)$ are depicted in Fig.~\ref{dydx050a}. It is seen that the transient time is about $30$~s, which conforms with the chosen value of decay parameter $\rho$ in \eqref{mk}, $\rho=\exp(-0.1T_s) = 0.9987$. Time histories of the transmission error $\delta_y$ and the coder range $M[k]$ are plotted in Fig.~\ref{delyMa} showing that the transmission error actually does not exceed $M[k]$. Time histories of $y_1(t)$ and $\bar {y}_1(t)$ for, $t\in[48, 50]$~s are depicted in Fig.~\ref{y1y1050a} to provide better imagination about the coding procedure. 

Normalized transmission error $Q_y$  as a function of the transmission rate $ R $ is plotted in Fig.~\ref{QyotRa}. Based on the simulation results the least-square estimate of the factor $G_y$ of the relation between the transmission error $Q_y$ and the channel rate $R$ having the inversely proportional form $Q_y=G_y/R$ was numerically found as $G_y=12.2$
 and the curve for $ G_y/R $ is also plotted in Fig.~\ref{QyotRa} (dashed line). 

Synchronization performance may be evaluated based on the normalized state synchronization error $Q$, \eqref{Q}, shown in Fig.~\ref{QotRa} as a function of the transmission rate $ R $. The simulation results make it possible to evaluate the parameter $G_y$ in the inversely proportional function $Q=G/R$. For the considered example $G=4.0$. The corresponding curve is plotted in Fig.~\ref{QotRa} (dashed line). One may notice that the synchronization error is less than the transmission error, $G<G_y$ (also compare Figs.~\ref{QyotRa} and \ref{QotRa}). This demonstrates the filtering abilities of the synchronization scheme. Based on a theoretical bound, simulation data are smoothened with a hyperbolic curve.
As seen from Fig.~\ref{QotKa} it is not necessary to increase control power (controller gain) in order to approach minimum of limit synchronization error. The value of the gain $K$ does not influence synchronization error for sufficiently large $K$.
\section{Conclusions}
Limit possibilities of controlled synchronization systems under
information constraints imposed by limited information capacity of
the coupling channel are evaluated. It is shown that the framework proposed in \cite{FradkovAndrievskyEvans_PRE06}, is suitable not only for observer-based synchronization but also for controlled master-slave synchronization via communication channel with limited information capacity.

 We propose a simple first order coder-decoder scheme and
 provide theoretical analysis for
multi-dimensional master-slave systems represented in the  Lurie
form. 
An output feedback control law is proposed based on the Passification theorem \cite{Fradkov74,FradkovMN99}. It is shown that upper and lower bounds
on the limit synchronization error are proportional to the maximum
rate of the coupling signal and inversely proportional to the
information transmission rate (channel capacity).  More complex coding procedures may provide better synchronization properties. 

The results are
applied to controlled synchronization of two chaotic Chua systems
coupled via a channel with limited capacity. Simulation results confirm theoretical investigations. 

Future research is aimed at examination of more complex system configurations, where control signal is subjected to information constraints too.

\section*{ACKNOWLEGEMENTS}
The work was supported by NICTA Victoria Research Laboratory of the University of Melbourne and the Russian Foundation for Basic
Research (proj. 05-01-00869, 06-08-01386).
 \begin{appendix}
\section*{Appendix}
Consider a linear system
\begin{align}
\dot{e}=Ae+B\xi(t), \quad\eps=Ce 
\label{A1}
\end{align}
with transfer function $W(\lambda)=C(\lambda I-A)^{-1}B$$=b(\lambda)/a(\lambda)$,
where $b(\lambda)$, $a(\lambda)$ are polynomials, degree of $a(\lambda)$ is $n$, degree of  $b(\lambda)$ is not greater than $n-1$. The system is called {\it hyper-minimum phase} (HMP), if $ b(\lambda)$ is  Hurwitz polynomial of degree  $n-1$ with positive coefficients. To find existence conditions for quadratic Lyapunov function we need the following result.

{\it Passification theorem \cite{Fradkov74,FradkovMN99}}. There exist positive-definite matrix $P=P\trn>0$ and a number $K$ such that
\begin{align}
PA_K+A_K\trn P<0,~~ PB=C\trn,~~ A_K=A-BKC  
\label{A1a}
\end{align}
if and only if   $W(\lambda)$ is HMP.

Consider a linear system with feedback
\begin{align}
\dot{e}=A_Ke+B\xi(t),~~ \eps=Ce,~~ A_K=A-BKC.
\label{A2}
\end{align}
Let us show that there exist a quadratic form $V(e)=e\trn Pe$ and a number $K$ such that time derivative $\dot V(e)$ of $V(e)$ along trajectories of  \eqref{A2} satisfies relation
\begin{align}
&\dot V(e)<0 ~~\text{for}~~ \xi\eps\ge 0, ~x\ne 0
\label{A3}
\end{align}
if and only if $W(\lambda)$ is HMP. To this end assume that $K$ is fixed. Relation \eqref{A3} is equivalent to existence of the matrix $P=P\trn>0$ such that $e\trn P(A_Ke+B\xi)+\xi Ce<0$ for $x\ne 0$. Since $\xi$ is arbitrary, the latter in turn, is equivalent to matrix relations  $PA_K+A_K\trn P<0$, $PB=C\trn$ and, by Passification theorem, to HMP condition. 

{\it Remark}. It also follows from Passification theorem that if HMP condition holds then $K$ satisfying \eqref{A1a} can be chosen sufficiently large. Besides, zero matrix in the right hand side of  the inequality in \eqref{A1a} can be replaced by matrix $-\mu P$ for sufficiently small $\mu>0$.

\end{appendix}

\end{document}